\documentclass[11pt]{article}

\usepackage[margin=1in]{geometry}
\usepackage{amsmath, amsthm, amssymb}
\usepackage[shortlabels]{enumitem}
\usepackage[numbers,sort&compress]{natbib}
\usepackage[colorlinks=true,allcolors=blue]{hyperref}

\theoremstyle{plain}
\newtheorem{theorem}{Theorem}[section]
\newtheorem{lemma}[theorem]{Lemma}
\newtheorem{proposition}[theorem]{Proposition}
\newtheorem{corollary}[theorem]{Corollary}
\theoremstyle{definition}
\newtheorem{remark}[theorem]{Remark}

\newtheorem{problem}{Problem}
\newtheorem{question}{Question}

\newcommand{\cP}{\mathcal{P}}
\newcommand{\BIP}{\mathrm{BIP}}
\newcommand{\PERF}{\mathrm{PERF}}
\newcommand{\GSP}{\mathrm{GSP}}
\newcommand{\UNIP}{\mathrm{UNIP}}
\newcommand{\coUNIP}{\mathrm{coUNIP}}
\newcommand{\cbip}{c_{\mathrm{BIP}}}
\newcommand{\cperf}{c_{\mathrm{PERF}}}
\newcommand{\cgsp}{c_{\mathrm{GSP}}}
\newcommand{\ccou}{c_{\mathrm{coUNIP}}}
\newcommand{\cuni}{c_{\mathrm{UNIP}}}
\newcommand{\ceil}[1]{\left\lceil #1\right\rceil}
\newcommand{\NN}{\mathbb{N}}
\DeclareMathOperator{\eq}{eq}

\newcommand{\doi}[1]{doi:\,\href{https://doi.org/#1}{\nolinkurl{#1}}}

\title
{Cover numbers by generalized split graphs and their subclasses}

\author{Donát Ákos Köller\thanks{Budapest University of Technology and Economics, \href{mailto:kollerda@math.bme.hu}{kollerda@math.bme.hu}} \and M\'{a}rton Marits\thanks{Yokohama National University, \href{mailto:marton-marits-xg@ynu.jp}{marton-marits-xg@ynu.jp}}}

\date{\today}
\usepackage[yyyymmdd]{datetime}

\begin{document}

\maketitle




\begin{abstract}
The cover number of a graph $G$ by a graph class $\cP$ is the minimum number of subgraphs of $G$ belonging to $\cP$ whose edge sets together cover the edge set of $G$. Marits (2026) asked for the exact values of the cover numbers by perfect, generalized split, co-unipolar and unipolar graphs. We show that the cover numbers by perfect, generalized split and co-unipolar graphs are not functions of the chromatic number $\chi$ and the clique number $\omega$. For all $n\ge w\ge2$ they attain the largest possible value $\lceil\log_2 n\rceil$ on a graph with $\chi=n$ and $\omega=w$, and for $\chi=\omega=n$ they attain every value from $1$ to $\lceil\log_2 n\rceil$. On the other hand, the cover numbers by co-unipolar and by generalized split graphs are at least $1+\lceil\log_2(\chi/\omega)\rceil$ and $\lceil\log_2(\chi/\omega)\rceil$, respectively, and these bounds are best possible up to small additive constants. For unipolar graphs we prove an exact formula: the cover number is the minimum, over all vertex sets $S$, of the larger of two quantities, namely the clique cover number of the subgraph induced by $S$ and the equivalence covering number of the subgraph induced by the remaining vertices. As consequences, we obtain an explicit formula for bipartite graphs, the exact values for complete bipartite graphs and for all hypercubes, a lower bound that is linear in $\chi$, and the NP-completeness of deciding whether three unipolar graphs suffice, even for triangle-free graphs.
\end{abstract}

\section{Introduction}\label{sec:intro}

All graphs in this paper are finite and simple, and all graph classes are closed under isomorphism. Let $\cP$ be a class of graphs. A \emph{covering} of a graph $G$ by $\cP$ is a family $G_1,\dots,G_k$ of subgraphs of $G$ such that $G_i\in\cP$ for each $i$ and $E(G)=E(G_1)\cup\dots\cup E(G_k)$. The \emph{cover number} $c_{\cP}(G)$ is the smallest $k$ for which such a covering exists. We use the convention that $c_{\cP}(G)=0$ if and only if $E(G)=\emptyset$. All graph classes considered in this paper contain 
$K_2$, so every graph has a covering by each of them and $c_{\cP}(G)\le|E(G)|$ is always finite.

Cover numbers include several classical parameters, for example arboricity (forests), thickness (planar graphs), and the equivalence covering number $\eq$, where the class consists of the disjoint unions of complete graphs \cite{Duchet1979,Alon1986,BlokhuisKloks1995,EsperetGimbelKing2010}. For surveys we refer to \cite{Harary1970,Pyber1991}.

Harary, Hsu and Miller \cite{HararyHsuMiller1977} showed that the cover number by bipartite graphs is $\cbip(G)=\ceil{\log\chi(G)}$, here and throughout, $\log$ denotes the binary logarithm unless stated otherwise. Coverings by comparability graphs were studied by Gy\'arf\'as, Marits and T\'oth \cite{GyarfasMaritsToth2026}.

Marits \cite{Marits2026} studied cover numbers by graph classes related to $\chi$-boundedness \cite{Gyarfas1987,ScottSeymour2020}. For every non-decreasing $f:\NN\to\NN$ with $f(n)\ge n$, he proved that
\begin{equation}\label{eq:marits3}
c_{\{\chi\le f(\omega)\}}(G)=\left\lceil{\frac{\log\chi(G)}{\log f(\omega(G))}} \right\rceil,
\end{equation}
where $\{\chi\le f(\omega)\}$ denotes the class of graphs $G$ with $\chi(G)\le f(\omega(G))$. He also established the chain of inequalities
\begin{equation}\label{eq:chain}
\left\lceil \frac{\log\chi}{\log\omega} \right\rceil=c_{\{\chi=\omega\}}\le\cperf\le\cgsp\le\ccou\le\cbip=\ceil{\log\chi},
\end{equation}
where $\PERF$, $\GSP$ and $\coUNIP$ denote the classes of perfect, generalized split and co-unipolar graphs (see Section~\ref{sec:prelim}). He proved that every gap in \eqref{eq:chain} can be arbitrarily large, and that the cover number $\cuni$ by unipolar graphs is not bounded from above by any function of $\chi$ and $\omega$. Among others, he posed the following two questions.

\begin{problem}[\cite{Marits2026}]\label{prob:1}
Determine the exact value of $\cperf$, $\ccou$ and $\cgsp$ in terms of the chromatic number or otherwise.
\end{problem}

\begin{problem}[\cite{Marits2026}]\label{prob:2}
Determine the exact value of $\cuni$.
\end{problem}

We provide the following solutions to these problems.

The first part of Problem~\ref{prob:1} asking for a formula in terms of the chromatic number is answered negatively, in the strongest possible sense. For each $\cP\in\{\PERF,\GSP,\coUNIP\}$, the upper bound $\ceil{\log\chi}$ in \eqref{eq:chain} is attained for every prescribed pair $(\chi,\omega)$, and when $\chi=\omega=n$, every value between $1$ and $\ceil{\log n}$ can be attained (Theorem~\ref{thm:nofunction}). Consequently, none of the three cover numbers can be expressed as a function of $\chi$ and $\omega$.

We then prove two lower bounds, on the cover number by co-unipolar graphs and by generalized split graphs. We have \[
    \ccou(G)\ge 1 + \left\lceil \log \left( \frac{\chi(G)}{\omega(G)}\right) \right\rceil
\] and \[
    \cgsp(G)\ge \left\lceil \log \left( \frac{\chi(G)}{\omega(G)}\right) \right\rceil,
\] proved in Theorems~\ref{thm:coUNIP} and~\ref{thm:GSP} respectively. Both are best possible up to small additive constants as functions of $\chi$ and $\omega$ (Proposition~\ref{prop:joins}). They determine $\ccou(G)$ up to at most $\ceil{\log\omega(G)}$ consecutive values, so for $\omega(G)\le4$, the value of $\ccou(G)$ is one of two consecutive integers determined by $\chi(G)$ (Corollary~\ref{cor:window}).

As for Problem \ref{prob:2}, we prove the exact formula \[
    \cuni(G)=\min_{S\subseteq V(G)}\max\{\theta(G[S]),\eq(G-S)\},
\] where $\theta$ is the clique cover number and $\eq$ is the equivalence covering number (Theorem~\ref{thm:unipformula}). As a corollary, we obtain an explicit formula involving the matching number for bipartite graphs (Corollary~\ref{cor:trianglefree}), the formula $\cuni(K_{m,n})=\min\{m,\ceil{n/2}\}$ for complete bipartite graphs with $m\le n$ (Proposition~\ref{prop:Kmn}), and the exact value of $\cuni(Q_d)$ for every hypercube, which extends Theorem~11 of \cite{Marits2026} from $d\ge 8$ to all $d$ (Theorem~\ref{thm:hypercube}). We also obtain a lower bound on $\cuni$ that is linear in $\chi$ (Proposition~\ref{prop:unipchi}), and the identity $\cuni(m\cdot H)=\eq(H)$ for $m\ge\eq(H)$ (Proposition~\ref{prop:copies}). From the latter we deduce that deciding whether $\cuni(G)\le 3$ is NP-complete, even for triangle-free graphs, whereas $\cuni(G)\le1$ can be decided in polynomial time (Theorem~\ref{thm:complexity}).

The paper is organized as follows. Section~\ref{sec:prelim} collects notation, definitions and some elementary facts that are used throughout. Sections~\ref{sec:p1} and~\ref{sec:p2} treat Problems~\ref{prob:1} and~\ref{prob:2}, respectively. Section~\ref{sec:open} lists open questions.

\section{Preliminaries}\label{sec:prelim}

For a graph $G$, we write $\chi(G)$, $\omega(G)$, $\alpha(G)$, $\Delta(G)$, $\chi'(G)$ and $\nu(G)$ for its chromatic number, clique number, independence number, maximum degree, chromatic index and matching number respectively. For unstated definitions, we refer the reader to the textbook \cite{diestel2025graph}. For a vertex $x$, $N(x)$ denotes its open neighborhood, and $N[x]=N(x)\cup\{x\}$. A set $D\subseteq V(G)$ \emph{dominates} $G$ if every vertex of $G$ lies in $D$ or has a neighbor in $D$, that is, $\bigcup_{x\in D}N[x]=V(G)$.

The \emph{clique cover number} $\theta(G)$ is the minimum number of cliques covering $V(G)$, equivalently, $\theta(G)=\chi(\overline{G})$. An \emph{equivalence graph} is a vertex-disjoint union of complete graphs, and the \emph{equivalence covering number} $\eq(G)$ is the minimum number of equivalence subgraphs of $G$ covering $E(G)$. See \cite{Duchet1979,Alon1986} for further discussion of this parameter. 

The sets of a partition of a set $V$ are called its \emph{cells}; we do not allow empty cells
unless stated otherwise. If $\pi_1,\dots,\pi_k$ are partitions of $V$, their
\emph{common refinement} is the partition of $V$ whose cells are the nonempty sets of the form
$C_1\cap\dots\cap C_k$ with $C_i\in\pi_i$. Equivalently, two elements of $V$ lie in the same
cell of the common refinement if and only if they lie in the same cell of every $\pi_i$.
For $k=0$, the common refinement is $\{V\}$. If $\pi$ is a partition of $V$ and
$U\subseteq V$, the \emph{restriction} of $\pi$ to $U$ is $\{C\cap U: C\in\pi\}$, with empty
sets removed.

We write $G\sqcup H$ for the disjoint union, $G\vee H$ for the join (the
disjoint union of $G$ and $H$ together with all edges between $V(G)$ and
$V(H)$), and $m\cdot H$ for the disjoint union of $m$ copies of $H$. The subgraphs $G_1,\dots,G_k$ of $G$ are said to \emph{cover} $G$ if $E(G)=\bigcup_i E(G_i)$.

A graph is \emph{unipolar} if its vertex set can be partitioned into sets $A$ and $B$ such that $A$ is a clique (the \emph{central clique}) and $G[B]$ is a disjoint union of complete graphs (the \emph{side cliques}). A graph is \emph{co-unipolar} if its complement is unipolar. Equivalently, $G$ is co-unipolar iff $V(G)=A\cup B$ where $A$ is independent and $B$ is partitioned into $p\ge0$ nonempty independent \emph{parts} $P_1,\dots,P_p$ such that any two vertices in different parts are adjacent; the edges between $A$ and $B$ are arbitrary. A graph is a \emph{generalized split graph} if it is unipolar or co-unipolar. We write $\BIP,\PERF,\UNIP,\coUNIP,\GSP$ for the classes of bipartite, perfect, unipolar, co-unipolar and generalized split graphs respectively.

All five of these classes are hereditary and closed under adding isolated vertices. Generalized split graphs are perfect \cite{PromelSteger1992,EschenWang2014}. Every bipartite graph with bipartition $(X,Y)$ is co-unipolar: take $A=X$ and $B=Y$ as a single part. Every complete graph is both unipolar and co-unipolar. Hence

\begin{equation}\label{eq:inclusions}
\BIP\subseteq\coUNIP\subseteq\GSP\subseteq\PERF,\qquad\text{so}\qquad \cperf\le\cgsp\le\ccou\le\cbip.
\end{equation}

Throughout, when $H$ is a unipolar graph, $A$ denotes its central clique and $\mathcal{Q}$ the family of its side cliques. When $H$ is co-unipolar, $A$ denotes the independent set and $P_1,\dots,P_p$ the parts, as in the definition above.

In a covering, a piece may be enlarged to a largest subgraph of $G$ with the same structure. For unipolar pieces, let $A$ be a clique of $G$ and $\mathcal{Q}$ a family of pairwise disjoint cliques of $G-A$, and let $U_{A,\mathcal{Q}}$ be the spanning subgraph of $G$ whose edges are the edges of $G[A]$, the edges of $G[Q]$ for $Q\in\mathcal{Q}$, and all edges of $G$ between $A$ and $V(G)\setminus A$. It is unipolar with central clique $A$ and side cliques $\mathcal{Q}$, and it contains every unipolar subgraph of $G$ with central clique $A$ and side cliques $\mathcal{Q}$.

For co-unipolar pieces, we call a partition $\pi=\{A,P_1,\dots,P_p\}$ of $V(G)$ \emph{admissible} if $p\ge0$, the cell $A$ may be empty, every $P_j$ is nonempty, and any two vertices lying in different sets $P_j\ne P_k$ are adjacent in $G$. Let $H_\pi$ be the spanning subgraph of $G$ whose edges are the edges of $G$ joining different cells of $\pi$. No edge of $H_\pi$ lies inside a cell, and the parts are pairwise completely joined, so $H_\pi$ is co-unipolar. Conversely, let $H$ be a co-unipolar subgraph of $G$ with $A$ and $P_1,\dots,P_p$. The parts are pairwise completely joined in $H$, hence in $G$, so $\pi=\{A\cup(V(G)\setminus V(H)),P_1,\dots,P_p\}$ is admissible. Every edge of $H$ joins $A$ to a part or joins two different parts, so $E(H)\subseteq E(H_\pi)$. Consequently,
\begin{equation}\label{eq:counipcells}
\ccou(G)\le k\iff E(G)\subseteq E(H_{\pi_1})\cup\dots\cup E(H_{\pi_k})\ \text{for some admissible partitions}\ \pi_1,\dots,\pi_k.
\end{equation}

\begin{lemma}\label{lem:basic}
Let $\cP$ be a hereditary class containing all stars. Then $c_{\cP}(G[S])\le c_{\cP}(G)$ and $c_{\cP}(G)\le c_{\cP}(G-v)+1$ for every graph $G$, every $S\subseteq V(G)$ and every $v\in V(G)$.
\end{lemma}

\begin{proof}
Let $G_1,\dots,G_k$ cover $G$. The graphs $G_i[S\cap V(G_i)]$ belong to $\cP$ and cover $G[S]$, because $E(G_i[S\cap V(G_i)])=E(G_i)\cap E(G[S])$. For the second inequality, add the star of $G$ at $v$ to a covering of $G-v$.
\end{proof}

Note that all classes in \eqref{eq:inclusions}, as well as $\UNIP$, contain all stars.

\begin{lemma}\label{lem:trianglefree}
If $Z$ is triangle-free, then every perfect subgraph of $Z$ is bipartite. Consequently $c_{\cP}(Z)=\ceil{\log\chi(Z)}$ for every $\cP\in\{\PERF,\GSP,\coUNIP\}$.
\end{lemma}

\begin{proof}
A perfect triangle-free graph $H$ has $\chi(H)=\omega(H)\le2$. Hence every covering of $Z$ by $\cP$ is a covering by bipartite graphs, so $c_\cP(Z)\ge\cbip(Z)$. The reverse inequality follows from \eqref{eq:inclusions}, and the value of $\cbip$ is given in \cite{HararyHsuMiller1977}.
\end{proof}

\section{Perfect, generalized split and co-unipolar covers}\label{sec:p1}

\begin{theorem}\label{thm:nofunction}
Let $\cP\in\{\PERF,\GSP,\coUNIP\}$ and let $2\le w\le n$ be integers.
\begin{enumerate}[(a)]
\item There is a graph $G$ with $\chi(G)=n$, $\omega(G)=w$ and $c_{\cP}(G)=\ceil{\log n}=\cbip(G)$.
\item For every integer $j$ with $1\le j\le\ceil{\log n}$ there is a graph $G$ with $\chi(G)=\omega(G)=n$ and $c_\cP(G)=j$.
\end{enumerate}
In particular, none of $\cperf$, $\cgsp$, $\ccou$ is a function of $\chi$ and $\omega$.
\end{theorem}

\begin{proof}
Let $Z$ be a triangle-free graph with $\chi(Z)=n$; for instance a Mycielski graph \cite{Mycielski1955}, or $K_2$ if $n=2$.

(a) Let $G=K_w\sqcup Z$. Then $\chi(G)=n$ and $\omega(G)=w$, since $\omega(Z)\le2\le w$. By Lemmas~\ref{lem:basic} and~\ref{lem:trianglefree}, $c_\cP(G)\ge c_\cP(Z)=\ceil{\log n}$. On the other hand, $c_\cP(G)\le\cbip(G)=\ceil{\log n}$ by \eqref{eq:inclusions}.

(b) Let $z_1,\dots,z_N$ be the vertices of $Z$, and let $G_i=(K_n\sqcup Z)-\{z_1,\dots,z_i\}$ for $0\le i\le N$. Every $G_i$ satisfies $\chi(G_i)=\omega(G_i)=n$. By (a) with $w=n$, $c_\cP(G_0)=\ceil{\log n}$. Since $K_n$ is co-unipolar, $c_\cP(G_N)=c_\cP(K_n)=1$. By Lemma~\ref{lem:basic}, consecutive values $c_\cP(G_i)$ and $c_\cP(G_{i+1})$ differ by at most one, so every integer between $1$ and $\ceil{\log n}$ occurs.
\end{proof}

Part (a) shows that the upper bound in \eqref{eq:chain} is attained on every level set $\{\chi=n,\ \omega=w\}$. Lemma~\ref{lem:trianglefree}, together with the theorem of Maffray and Preissmann \cite{MaffrayPreissmann1996} that $4$-colorability of triangle-free graphs is NP-complete, also shows that deciding $c_\cP(G)\le 2$ is NP-hard for $\cP\in\{\PERF,\GSP,\coUNIP\}$, even for triangle-free $G$.

The lower bound $\ceil{\log\chi/\log\omega}$ in \eqref{eq:chain} comes from the inclusion $\GSP\subseteq\{\chi=\omega\}$. For co-unipolar and generalized split graphs, the structure of the pieces gives much stronger bounds when $\chi$ is large compared with $\omega$.

\begin{theorem}\label{thm:coUNIP}
Let $t \geq 1$. Every graph $G$ with $\ccou(G) \leq t$ satisfies $\chi(G)\le\omega(G)\,2^{t-1}$. Consequently, every graph $G$ with at least one edge satisfies \[
    \ccou(G) \geq 1 + \left\lceil \log \left( \frac{\chi(G)}{\omega(G)} \right) \right\rceil.
\]
\end{theorem}

\begin{proof}
We use induction on $t$. If $t=1$, then $G$ is co-unipolar, hence it is also perfect, which implies that $\chi(G)=\omega(G)$.

Let $t\ge2$. If $\ccou(G)\le t-1$ we are done by induction, so assume that co-unipolar subgraphs $H_1,\dots, H_t$ cover $G$. As shown before \eqref{eq:counipcells}, there is an admissible partition $\pi=\{A,P_1,\dots,P_p\}$ with $E(H_1)\subseteq E(H_\pi)$. Let $C$ be a nonempty cell of $\pi$. Since $E(H_1)\subseteq E(H_\pi)$, every edge of $H_1$ joins two different cells of $\pi$; an edge of $G[C]$ has both ends in $C$, so it is not an edge of $H_1$. Hence $G[C]$ is covered by the graphs $H_i[C\cap V(H_i)]$ with $2\le i\le t$, which are co-unipolar since $\coUNIP$ is hereditary. So $\ccou(G[C])\le t-1$, and by induction $\chi(G[C])\le\omega(G[C])\,2^{t-2}$.

Since $\pi$ is admissible, any two vertices lying in different parts $P_j\ne P_k$ are adjacent in $G$. Hence, if $K_j$ is a maximum clique of $G[P_j]$ for $j=1,\dots,p$, then $K_1\cup\dots\cup K_p$ is a clique of $G$, which gives $\sum_j\omega(G[P_j])\le\omega(G)$. Coloring the cells with pairwise disjoint palettes, we obtain
\begin{equation}\label{eq:cellsum}
\chi(G)\le\chi(G[A])+\sum_{j=1}^p\chi(G[P_j])\le2^{t-2}\Bigl(\omega(G)+\sum_{j=1}^p\omega(G[P_j])\Bigr)\le\omega(G)\,2^{t-1},
\end{equation}
where the term for $A$ is omitted if $A=\emptyset$.

For the consequence, let $t=\ccou(G)\ge1$. Then $t-1\ge\log(\chi(G)/\omega(G))$, and since $t-1$ is an integer, $t-1\ge\ceil{\log(\chi(G)/\omega(G))}$.
\end{proof}

The bound of Theorem~\ref{thm:coUNIP} is attained when $\omega=2$, and it is sharp up to a factor of $2$ in $\chi$ (equivalently, an additive $1$ in $\ccou$) for every even $\omega$.

\begin{proposition}\label{prop:joins}
Let $r,s\ge1$, let $Z_1,\dots,Z_r$ be triangle-free graphs with $\chi(Z_i)=2^s$, and let $J=Z_1\vee\dots\vee Z_r$. Then $\omega(J)=2r$, $\chi(J)=r2^s$ and $s\le\cgsp(J)\le\ccou(J)\le s+1$. In particular $\chi(J)\ge\omega(J)\,2^{\ccou(J)-2}$, and for $r=1$ we have $\chi(J)=\omega(J)\,2^{\ccou(J)-1}$.
\end{proposition}

\begin{proof}
The chromatic and clique numbers are additive under joins, and $\omega(Z_i)=2$. The lower bound follows from Lemmas~\ref{lem:basic} and~\ref{lem:trianglefree} applied to the induced subgraph $Z_1$.

For the upper bound, $\pi_0=\{\emptyset,V(Z_1),\dots,V(Z_r)\}$ is admissible, because the sets $V(Z_i)$ are pairwise completely joined. Fix proper colorings of the $Z_i$ with colors in $\{0,1\}^s$; together they define a map $\varphi:V(J)\to\{0,1\}^s$. For $\ell=1,\dots,s$, let $\pi_\ell$ consist of the cell $\{v:\varphi(v)_\ell=0\}$ and the single part $\{v:\varphi(v)_\ell=1\}$; this partition is trivially admissible. An edge of $J$ between different $Z_i$ joins different cells of $\pi_0$. An edge inside some $Z_i$ joins vertices whose colors differ in some coordinate $\ell$, so it joins different cells of $\pi_\ell$. By \eqref{eq:counipcells}, $\ccou(J)\le s+1$.

For $r=1$ we have $J=Z_1$, so $\ccou(J)=s$ by Lemma~\ref{lem:trianglefree}, and $\chi(J)=2^s=\omega(J)\,2^{s-1}$.
\end{proof}

A covering by generalized split graphs may mix unipolar and co-unipolar pieces. The central cliques of the unipolar pieces are few and small, and outside them the unipolar pieces only contribute bounded degree. This leads to the following bound.

\begin{theorem}\label{thm:GSP}
Let $t\ge1$. Every graph $G$ with $\cgsp(G)\le t$ satisfies $\chi(G)\le\omega(G)\,(2^{t-1}+t)\le\omega(G)\,2^{t}$. Consequently, \[
    \cgsp(G) \geq \left\lceil \log\left(\frac{\chi(G)}{\omega(G)}\right) \right\rceil
\] for every graph $G$.
\end{theorem}

\begin{proof}
If $G$ has no edges there is nothing to prove. Otherwise let the unipolar graphs $U_1, \dots, U_u$ and the co-unipolar graphs $H_1, \dots, H_k$ cover $G$. Since $2^{t-1}+t$ is increasing in $t$, we can assume that $u + k = t = \cgsp(G)$. If $u=0$, Theorem~\ref{thm:coUNIP} gives $\chi(G)\le\omega(G)2^{t-1}$. So assume $u \geq 1$ (and therefore $t \geq 1$), and write $\omega := \omega(G)$.

Let $A_i$ and $\mathcal{Q}_i$ be the central clique and the family of side cliques of $U_i$. As shown before \eqref{eq:counipcells}, there are admissible partitions $\pi_1,\dots,\pi_k$ with $E(H_j)\subseteq E(H_{\pi_j})$. Let $S=A_1\cup\dots\cup A_u$ and $G'=G-S$. Let $\pi'_j$ be the restriction of $\pi_j$ to $V(G')$, which is an admissible partition of $G'$.

\emph{Claim 1.} Let $C$ be a cell of the common refinement of $\pi'_1,\dots,\pi'_k$. Then $\Delta(G[C])\le u(\omega(G[C])-1)$, hence $\chi(G[C])\le u\,\omega(G[C])$.

Let $vw$ be an edge of $G[C]$. Since $v,w$ lie in the same cell of every $\pi_j$, the edge $vw$ is not in any $H_{\pi_j}$, so it belongs to some $U_i$. As $v,w\notin S$, the edge $vw$ is not incident with $A_i$, so it lies inside a side clique $Q\in\mathcal{Q}_i$. The vertex $v$ lies in at most one side clique $Q_i(v)$ of each $U_i$. Therefore every neighbor of $v$ in $G[C]$ lies in some $(Q_i(v)\cap C)\setminus\{v\}$. Each $Q_i(v)\cap C$ is a clique of $G[C]$, so $v$ has at most $u(\omega(G[C])-1)$ neighbors in $G[C]$. A greedy coloring then uses at most $u(\omega(G[C])-1)+1\le u\,\omega(G[C])$ colors.

\emph{Claim 2.} Let $F$ be a graph and $\rho_1,\dots,\rho_\kappa$ admissible partitions of $V(F)$. Suppose $\chi(F[C])\le u\,\omega(F[C])$ for every cell $C$ of their common refinement. Then $\chi(F)\le u\,\omega(F)\,2^\kappa$.

We use induction on $\kappa$; the case $\kappa=0$ is the hypothesis. Let $\rho_1=\{A,P_1,\dots,P_p\}$. The restrictions of $\rho_2,\dots,\rho_\kappa$ to a cell of $\rho_1$ are admissible, and their common refinement consists of cells of the full refinement. By induction, as in \eqref{eq:cellsum},
\begin{equation}\label{eq:claim2}
    \chi(F)\le\chi(F[A])+\sum_j\chi(F[P_j])\le u2^{\kappa-1}\Bigl(\omega(F)+\sum_j\omega(F[P_j])\Bigr)\le u\,\omega(F)\,2^\kappa.
\end{equation}

By Claims 1 and 2, $\chi(G')\le u\,\omega\,2^{k}$. Moreover $\chi(G[S])\le|S|\le u\omega$. Hence $\chi(G)\le\chi(G[S])+\chi(G')\le u\omega(1+2^{t-u})$. Since $u\le2^{u-1}$ for $u\ge1$, we have $u2^{t-u}\le2^{t-1}$. Also $u\le t$. Therefore $\chi(G)\le\omega(2^{t-1}+t)$, which proves the theorem.
\end{proof}

For the graphs $J$ of Proposition~\ref{prop:joins}, built from $r$ triangle-free graphs with chromatic number $2^s$, we have $\chi(J)/\omega(J)=r2^s/(2r)=2^{s-1}$. Thus Theorem~\ref{thm:GSP} gives $\cgsp(J)\ge s-1$, while Proposition~\ref{prop:joins} gives $\cgsp(J)\le s+1=\log(\chi(J)/\omega(J))+2$. Hence the bound of Theorem~\ref{thm:GSP} cannot be improved by more than an additive constant~$2$.

\begin{corollary}\label{cor:window}
Every graph $G$ with at least one edge satisfies
\begin{equation}
\left\lceil \log \left( \frac{\chi(G)}{\omega(G)}\right) \right\rceil \leq \cgsp(G)\le\ccou(G)\le\ceil{\log\chi(G)},\label{eq:window}
\end{equation}
and
\begin{equation}
\ccou(G) \geq 1 + \left\lceil \log \left( \frac{\chi(G)}{\omega(G)}\right) \right\rceil.\label{eq:window2}
\end{equation}
In particular, given $\chi(G)$ and $\omega(G)$, the value $\ccou(G)$ lies in a set of at most $\ceil{\log\omega(G)}$ consecutive integers. Thus $\ccou(G)=\ceil{\log\chi(G)}$ if $\omega(G)=2$, and $\ccou(G)\in\{\ceil{\log\chi(G)}-1,\ceil{\log\chi(G)}\}$ if $\omega(G)\le4$.
\end{corollary}

\begin{proof}
Combine Theorems~\ref{thm:coUNIP} and~\ref{thm:GSP} with \eqref{eq:chain}. Writing $x=\log\chi(G)$ and $y=\log\omega(G)$, the number of possible values of $\ccou(G)$ is $\ceil{x}-\ceil{x-y}\le\ceil{y}$.
\end{proof}

\begin{remark}
Let $w\ge2$ and $t\ge2$, and let $G$ be a graph with $\omega(G)=w$ and $\chi(G)=w^t$ (for instance, the disjoint union of $K_w$ with a graph of clique number at most $w$ and chromatic number $w^t$). The lower bound in \eqref{eq:chain} gives $\ccou(G)\ge c_{\{\chi=\omega\}}(G)=t$. Theorem~\ref{thm:coUNIP} gives $\ccou(G)\ge1+(t-1)\log w$, which is larger by a factor of about $\log w$.
\end{remark} 

For perfect graphs, the only general bound of this type we are aware of is $\chi\le\omega^t$ for unions of $t$ perfect graphs, which yields the lower bound in \eqref{eq:chain}. Whether it can be replaced by a bound of the form $\chi\le C^t\omega$ is open (Section~\ref{sec:open}).

\section{Unipolar covers}\label{sec:p2}

Let us now move on to coverings by unipolar graphs. We have proved the following exact formula for the cover number by unipolar graphs.

\begin{theorem}\label{thm:unipformula}
For every graph $G$,
\begin{equation}\label{eq:unipformula}
\cuni(G)=\min_{S\subseteq V(G)}\max\bigl\{\theta(G[S]),\ \eq(G-S)\bigr\}.
\end{equation}
\end{theorem}

\begin{proof}
($\ge$) Let $U_1,\dots,U_t$ be unipolar subgraphs covering $G$, with central cliques $A_1,\dots,A_t$. These are cliques of $G$, so $S=A_1\cup\dots\cup A_t$ satisfies $\theta(G[S])\le t$. For each $i$, let $E_i$ be the subgraph of $U_i$ induced by $V(U_i)\setminus S$. This is an induced subgraph of $U_i-A_i$, which is a disjoint union of complete graphs, so $E_i$ is an equivalence subgraph of $G-S$. Every edge of $G-S$ lies in some $U_i$, and since both its ends are outside $S\supseteq A_i$, it lies in $E_i$. Hence $\eq(G-S)\le t$.

($\le$) Let $S\subseteq V(G)$ and $t=\max\{\theta(G[S]),\eq(G-S)\}$. Choose cliques $A_1,\dots,A_t$ (some possibly empty) covering $S$, and equivalence subgraphs $E_1,\dots,E_t$ of $G-S$ (some possibly edgeless) covering $G-S$. Let $\mathcal{Q}_i$ be the family of vertex sets of the components of $E_i$; these are cliques of $G-S$, and hence of $G-A_i$. Put $U_i=U_{A_i,\mathcal{Q}_i}$, which is unipolar.

Every edge with an end in $S$ has an end in some $A_i$, so it lies in $U_i$. Every edge of $G-S$ lies in some $E_i$, hence inside a member of $\mathcal{Q}_i$, hence in $U_i$. Thus $\cuni(G)\le t$.
\end{proof}

Taking $S=V(G)$ and $S=\emptyset$ in \eqref{eq:unipformula} gives
$\cuni(G)\le\min\{\theta(G),\eq(G)\}$. A proper edge coloring of $G$
with $\chi'(G)$ colors splits $E(G)$ into $\chi'(G)$ matchings, and every
matching is an equivalence graph, so $\eq(G)\le\chi'(G)$. Together with
Vizing's theorem \cite{diestel2025graph}, we obtain
\begin{equation}\label{eq:unipupper}
\cuni(G)\le\min\{\theta(G),\eq(G)\}\le\chi'(G)\le\Delta(G)+1.
\end{equation}

Both parameters and the minimum over $S$ are needed. Let $n\ge2$ and $G=K_{1,n}\sqcup n\cdot K_2$. Then $\theta(G)=2n$ and $\eq(G)=n$, but $\cuni(G)=1$: take $S$ to be the center of the star, or observe directly that $G$ is unipolar.

\begin{corollary}\label{cor:trianglefree}
    The following two statements hold.
    \begin{enumerate}
        \item If $G$ is triangle-free, then
            \begin{equation}\label{eq:tf}
                \cuni(G)=\min_{S\subseteq V(G)}\max\bigl\{|S|-\nu(G[S]),\ \chi'(G-S)\bigr\}.
            \end{equation}
        \item If $G$ is bipartite, then
            \begin{equation}\label{eq:bip}
                \cuni(G)=\min\bigl\{t\ge0:\ \exists S\subseteq V(G)\ \text{with}\ |S|-\nu(G[S])\le t\ \text{and}\ \Delta(G-S)\le t\bigr\}.
            \end{equation}
    \end{enumerate}
\end{corollary}

\begin{proof}
In a triangle-free graph $F$, a clique cover consists of edges and single vertices. 
In a cover that uses the fewest possible edges, the edges used necessarily form a matching, hence $\theta(F)=|V(F)|-\nu(F)$. The equivalence subgraphs of $F$ are matchings. A covering by matchings can be turned into a partition into matchings, so $\eq(F)=\chi'(F)$. For bipartite $F$, $\chi'(F)=\Delta(F)$ by K\H{o}nig's theorem \citep{diestel2025graph}.
\end{proof} 

In other words, a bipartite graph satisfies $\cuni(G)\le t$ if and only if there is a vertex set $S$, covered by at most $t$ edges and isolated vertices, whose deletion can reduce the maximum degree of $G$ to at most $t$.

\begin{proposition}\label{prop:Kmn}
For $1\le m\le n$, we have $\cuni(K_{m,n})=\min\{m,\ceil{n/2}\}$.
\end{proposition}




\begin{proof}
Let $X$ and $Y$ be the two sides of $G=K_{m,n}$, with $|X|=m$ and $|Y|=n$. For $S\subseteq V(G)$, write $a=|S\cap X|$ and $b=|S\cap Y|$. Then $G[S]=K_{a,b}$, so $|S|-\nu(G[S])=\max\{a,b\}$. We use \eqref{eq:bip}.

Taking $S=X$ gives the value $m$. Taking $a=\lfloor m/2\rfloor$ and $b=\lfloor n/2\rfloor$ gives $\max\{a,b\}=\lfloor n/2\rfloor$ and $G-S=K_{\ceil{m/2},\ceil{n/2}}$, whose maximum degree is $\ceil{n/2}$. So $\cuni(G)\le\min\{m,\ceil{n/2}\}$.

Conversely, let $S$ witness $\cuni(G)\le t$ in \eqref{eq:bip}. If $X\subseteq S$, then $t\ge a=m$. Otherwise, a vertex of $X$ outside $S$ has $n-b$ neighbors in $G-S$. Hence $t\ge\max\{b,n-b\}\ge n/2$, and so $t\ge\ceil{n/2}$, as $t$ is an integer.
\end{proof}


In particular $\cuni(K_{n,n})=\ceil{n/2}$, which gives a simple family with $\chi=\omega=2$ and unbounded unipolar cover number. For the hypercubes $Q_d$, Marits \cite[Theorem~11]{Marits2026} proved $\cuni(Q_d)=d$ for $d\ge8$. We determine the exact value of this cover number for all $d$.

\begin{theorem}\label{thm:hypercube}
    We have \begin{equation}
        \cuni(Q_d) = \begin{cases}
                d & \text{if } d = 1 \text{ or } d \geq 6 \\
                d- 1 & \text{if } d \in \{2,3,4,5\}
        \end{cases}
    \end{equation}
\end{theorem}

\begin{proof}
$Q_d$ is $d$-regular and bipartite on $2^d$ vertices, so $\cuni(Q_d)\le d$ by \eqref{eq:bip} with $S=\emptyset$. Suppose $S$ witnesses $\cuni(Q_d)\le t$ with $t<d$. Then $S$ is the vertex set of at most $t$ edges and vertices, so $|S|\le2t$. Moreover, a vertex $v\notin S$ has $d$ neighbors in $Q_d$, of which at most $\Delta(Q_d-S)\le t$ lie outside $S$; hence $v$ has at least $d-t\ge1$ neighbors in $S$. We count the edges between $S$ and $V(Q_d)\setminus S$ in two ways. Each of the $2^d-|S|$ vertices outside $S$ sends at least $d-t$ edges into $S$, and each vertex of $S$ sends at most $d$ edges out of $S$. Using also $|S|\le2t$ and $d-t>0$, we get
\begin{equation}\label{eq:qdcount}
(2^d-2t)(d-t)\le(2^d-|S|)(d-t)\le d|S|\le2td.
\end{equation}

\emph{The case $d\ge6$.} It suffices to exclude $t=d-1$. Then every vertex outside $S$ has a neighbor in $S$, so the closed neighborhoods $N[x]$, $x\in S$, cover $V(Q_d)$. Since $Q_d$ is triangle-free, $S$ is the union of at most $d-1$ sets, each a single vertex or the two ends of an edge. A single vertex $x$ has $|N[x]|=d+1\le2d$. For an edge $uv$, the sets $N[u]$ and $N[v]$ have $d+1$ elements each and share only $u$ and $v$, as a further common element would form a triangle; hence $|N[u]\cup N[v]|=2d$. Therefore $2d(d-1)\ge2^d$.

\emph{Lower bounds for $d\le5$.} Inequality \eqref{eq:qdcount} fails for $(d,t)=(3,1)$, $(4,2)$ and $(5,3)$, since $12>6$, $24>16$ and $52>30$. So $\cuni(Q_d)\ge d-1$ for $d\in\{3,4,5\}$.

\emph{Upper bounds for $d\le5$.} $Q_1=K_2$ and $Q_2=K_{2,2}$ are unipolar (Proposition~\ref{prop:Kmn}). For $Q_3$, take $S=\{000,111\}$: then $\theta(Q_3[S])=2$ and $Q_3-S$ is a $6$-cycle. For $Q_4$, take $S=\{0000,0001,1110,1111\}$, the vertex set of two edges. Vertices of weight $1$ and $3$ are adjacent to $0000$ and $1111$, respectively. A vertex of weight $2$ ending in $1$ is adjacent to $0001$, and one ending in $0$ is adjacent to $1110$. So $S$ dominates $Q_4$ and $\Delta(Q_4-S)\le3$. For $Q_5$, let $S$ be the set of the eight vertices $x0$ and $x1$ with $x$ in the previous set. It is the vertex set of four edges, and it dominates $Q_5$ because the previous set dominates $Q_4$. Hence $\Delta(Q_5-S)\le4$.
\end{proof}

Marits \cite{Marits2026} showed that $\cuni$ admits no upper bound in terms of $\chi$ and $\omega$. Theorem~\ref{thm:unipformula} yields a lower bound that is linear in $\chi$, in contrast with the logarithmic bounds of Section~\ref{sec:p1}.

\begin{proposition}\label{prop:unipchi}
Every graph $G$ with $\cuni(G)\le t$ satisfies $\chi(G)\le t(2\omega(G)-1)+1$, hence \[
    \cuni(G)\ge \frac{\chi(G)-1}{2\omega(G)-1}.
\]
\end{proposition}

\begin{proof}
By Theorem~\ref{thm:unipformula} there is a set $S\subseteq V(G)$ with $\theta(G[S])\le t$ and $\eq(G-S)\le t$. Write $\omega=\omega(G)$. Since $\theta(G[S])\le t$, we have $\chi(G[S])\le|S|\le t\omega$. On the other hand, in $G-S$, a vertex lies in at most one clique of each of the $t$ equivalence subgraphs, and each such clique has at most $\omega$ vertices. So $\Delta(G-S)\le t(\omega-1)$ and $\chi(G-S)\le t(\omega-1)+1$.

Coloring $G[S]$ and $G-S$ with disjoint sets of colors yields a proper coloring of $G$, therefore \begin{equation}\label{eq:unipchi} \chi(G)\le\chi(G[S])+\chi(G-S)\le t\omega+t(\omega-1)+1=t(2\omega-1)+1, \end{equation} which proves the proposition.
\end{proof}

\begin{proposition}\label{prop:copies}
For every graph $H$ and every integer $m\ge\eq(H)$, $\cuni(m\cdot H)=\eq(H)$.
\end{proposition}

\begin{proof}
Let $q=\eq(H)$; if $q=0$ there is nothing to prove. Taking the union of the copies of each of $q$ equivalence subgraphs covering $H$ gives $q$ equivalence subgraphs covering $m\cdot H$, so $\cuni(m\cdot H)\le\eq(m\cdot H)\le q$ by \eqref{eq:unipupper}.

Suppose that $t=\cuni(m\cdot H)<q$. By Theorem~\ref{thm:unipformula} there is a set $S\subseteq V(m\cdot H)$ with $\theta((m\cdot H)[S])\le t$ and $\eq((m\cdot H)-S)\le t$. Every clique of $m\cdot H$ lies within one copy of $H$, and $S$ is covered by at most $t<q\le m$ cliques, so some copy $H'$ of $H$ is disjoint from $S$. Then $H'$ is an induced subgraph of $(m\cdot H)-S$. Restricting equivalence subgraphs to an induced subgraph yields equivalence subgraphs, so $q=\eq(H)=\eq(H')\le\eq((m\cdot H)-S)\le t<q$, a contradiction.
\end{proof}

Since $\eq(H)\le|E(H)|$, Proposition~\ref{prop:copies} reduces the computation of $\eq$ to that of $\cuni$. Computing $\eq$ is NP-hard even for split graphs \cite{BlokhuisKloks1995}, so the same holds for $\cuni$. More precisely, we have the following.

\begin{theorem}\label{thm:complexity}
\leavevmode
\begin{enumerate}[(a)]
\item Deciding whether $\cuni(G)\le1$ is in P.
\item Deciding whether $\cuni(G)\le3$ is NP-complete, even for triangle-free graphs $G$.
\end{enumerate}
\end{theorem}

\begin{proof}
(a) Deciding whether $\cuni(G)\le1$ is equivalent to deciding whether $G \in \UNIP$, and it is known from \cite{EschenWang2014} and \cite{McDiarmidYolov2015} that this problem is decidable in polynomial time.

(b) A certificate consists of three pairs $(A_i,\mathcal{Q}_i)$, each a clique and a family of pairwise disjoint cliques of $G-A_i$; whether the graphs $U_{A_i,\mathcal{Q}_i}$ cover $G$ can be checked in polynomial time. For hardness, we use the result of Kami\'nski and Lozin \cite{KaminskiLozin2007} (extending Holyer \cite{Holyer1981}) that edge $3$-colorability is NP-complete for cubic graphs of girth at least $4$. Let $H$ be such a graph. By Corollary~\ref{cor:trianglefree} and Vizing's theorem, $\eq(H)=\chi'(H)\in\{3,4\}$. By Proposition~\ref{prop:copies} with $m=4$, the triangle-free graph $4\cdot H$ satisfies $\cuni(4\cdot H)=\chi'(H)$. Hence $\cuni(4\cdot H)\le3$ if and only if $H$ is $3$-edge-colorable.
\end{proof}

\section{Concluding remarks and open problems}\label{sec:open}

Our results settle the ``chromatic number'' part of Problem~\ref{prob:1} negatively and give bounds in terms of $\chi$ and $\omega$ for $\ccou$ and $\cgsp$ that are best possible up to additive constants. They also reduce Problem~\ref{prob:2} to the classical parameters $\theta$ and $\eq$. Because $\cuni(m\cdot H)=\eq(H)$, any exact formula for $\cuni$ must contain the equivalence covering number in some form, and Theorem~\ref{thm:unipformula} shows that nothing more is needed. 

\begin{question}
Is there a constant $C$ such that every graph covered by $t$ perfect graphs satisfies $\chi\le C^t\omega$? An affirmative answer would give $\cperf\ge\log(\chi/\omega)/\log C$.
\end{question}

\begin{question}
Determine $g(w,t)=\max\{\chi(G):\ \omega(G)\le w,\ \ccou(G)\le t\}$. Theorem~\ref{thm:coUNIP} and Proposition~\ref{prop:joins} give $w2^{t-2}\le g(w,t)\le w2^{t-1}$ for even $w$ and $t\ge2$, with equality on the right for $w=2$.
\end{question}

\begin{question}
What is the complexity of deciding $\cuni(G)\le2$? Is $\cuni$ polynomial-time computable on bipartite graphs, where \eqref{eq:bip} reduces it to a degree-and-matching problem? What is the best constant in Proposition~\ref{prop:unipchi} for $\omega\ge3$?
\end{question}

\section*{Acknowledgements}
The authors used Claude Opus 5.5 to develop and check the results in the paper. The authors independently reviewed and edited the outputs and take full responsibility for its content.

\bibliographystyle{plainnat}
\bibliography{refs}

\end{document}